# On transformation formulae for Srivastava-Daoust type *q*-hypergeometric series

**Yashoverdhan Vyas*, Kalpana Fatawat**

Department of Mathematics, School of Engineering, Sir Padampat Singhania University Udaipur, Rajasthan.

**Abstract:** We present here the *q*-analogues of certain transformations or reduction formulae for Srivastava-Daoust type double hypergeometric series. These reduction formulae are derived by utilizing the extended Bailey's Transform developed and studied by Joshi and Vyas [Int. J. Math. Sci., (12), 2005, 1909-1927]. A number of well-known *q*-hypergeometric transformations are also obtained as special cases of our results.

**Key words**: Extended Bailey's Transform, *q*-analogue, Hypergeometric Series, Srivastava-Daoust Series, Reduction Formulae.

## INTRODUCTION

The enormous success of the theory of hypergeometric series in single variable has stimulated the development of a corresponding theory in two or more variables. Prior to 1880 when Appell investigated double hypergeometric series, Hermite (1865) introduced some polynomials which are particular cases of Appell's double series $F_3$ but the credit of first systematic study of multiple hypergeometric series goes to Appell. The first Appell hypergeometric function $F_1$ of two variables is given by

$$F_1[\alpha;\beta,\beta';\gamma;u,v]=\sum_{m,n\geq 0}\frac{(\alpha)_{m+n}(\beta)_m(\beta')_n u^m v^n}{(\gamma)_{m+n}\, m!n!}$$

This Appell function $F_1$ is a solution of the system of partial differential equation given below.

$$u(1-u)\frac{\partial^2\Omega}{\partial u^2}+v(1-u)\frac{\partial^2\Omega}{\partial u\partial v}+[\gamma-(\alpha+\beta+1)u]\frac{\partial\Omega}{\partial u}-\beta v\frac{\partial\Omega}{\partial v}-\alpha\beta\Omega=0$$
$$v(1-v)\frac{\partial^2\Omega}{\partial v^2}+u(1-v)\frac{\partial^2\Omega}{\partial u\partial v}+[\gamma-(\alpha+\beta'+1)v]\frac{\partial\Omega}{\partial v}-\beta' u\frac{\partial\Omega}{\partial u}-\alpha\beta'\Omega=0 \qquad (1.1)$$

satisfying $\Omega=\Omega(u,v)=F_1(\alpha,\beta,\beta';\gamma;u,v)$

Besides this, the remaining three Appell hypergeometric functions $F_2, F_3$ and $F_4$ in two variables are also solutions of the system of partial differential equations which are recorded in Slater[1]. In 1893, Le Vavasseur presented a tableau of 60 convergent solutions of the system (1.1) in terms of Appell's double hypergeometric series $F_1$. The tableau was reproduced by Appell and Kampé de Fériet in the monograph [2] along with references to relevant literature on the subject and a summary of the important results concerning (1.1) obtained by eminent mathematicians e.g. Pochhammer, Picard, Goursat and others. Later on, Erdélyi presented a systematic investigation of contour integrals satisfying equation (1.1) and thereby obtained the fundamental set of solutions (including 60 solutions in terms of Horn's series $G_2$, which represents 15 new distinct solutions) for vicinity of all singular points of (1.1) The complete account of development of multiple hypergeometric series along with its applications is given in Srivastava and Karlsson [3] and Exton[4].

A large number of problems from science and engineering can be represented in terms of differential equations, and their solutions can be put in the form of multiple hypergeometric functions or their limiting cases. Hence, it is worthwhile to study the multiple hypergeometric series, their reductions or transformations and summations and other types of relations. Exton[4], Srivastava and Daoust [5-7], Qureshi *et al.*[8], Hai *et al.* [9], Bushman and Srivastava[10] have contributed in the field of multiple hypergeometric series and also discussed their applications.

It is well known that whenever a generalized hypergeometric function reduces to quotient of the products of the gamma function, the results are very important from the application point of view in numerous areas of physics, mathematics and statistics including in series systems of symbolic computer algebra manipulation. Similarly, the transformation or reduction formulas for certain classes of multiple (especially double) series are very useful in a number of applications. For example, Srivastava[11] investigated some of the reduction and summation formulas for generalized multiple hypergeometric series that arises naturally in physical and quantum chemical problems. Singh[12] evaluated three integrals involving Kampé de Fériet function and three expansion formula for Kampé de Fériet function. He found the application of these results in solving boundary value problem (heat equation) and in the derivation of radial wave functions respectively. Srivastava *et al.*[13,14] developed several reduction formulae for the double hypergeometric series .

Furthermore, many of the researchers have found the basic analogues of multiple hypergeometric series and their reduction or transformation formulae and also discussed their applications. For details, we refer, Andrews *et al.*[15], Andrews[16], Gasper and Rahman[17], Exton[4], Srivastava and Manocha[18], Srivastava and Karlsson[3], Saxena and Gupta[19], Ernst[20] and references therein.

Recently, the authors [21] investigated fourteen Srivastava-Daoust type reduction formulae using extended Bailey transform technique investigated by Joshi and Vyas[22]. In this paper, our aim is to develop the $q$-extensions of the reduction formulae obtained in[21]. Many of the derived reductions or transformations are interesting extensions of the $q$-analogues of some well-known results in the field of hypergeometric series e.g. Euler transformation formula, Whipple's quadratic transformation and one of the Kampé de Fériet reductions given in Srivastava and Karlsson[3]. We apply many of the well known $q$-hypergeometric identities and the $q$-Pfaff-Saalschütz summation theorem. For further details on $q$-hypergeometric identities and notation, we refer, Gasper and Rahman[19] and Ernst[20].

The generalized reduction formulae listed in section two (from (2.2) to (2.4)) are written explicitly without using any notation (except (2.1) and (2.3)), since the available notations of the $q$-Srivastava-Daoust series (1.5) are not appropriate to express such results containing additional powers of $q$. However, the single $q$-hypergeometric series given in equation (1.3) have notations for such results containing additional powers of $q$ and the notations for the $q$-Srivastava-Daoust series given in (1.5) can be developed on the line of the notations given for the Single $q$-hypergeometric series in (1.3) and the $q$-Kampé de Fériet series in (1.6). But these notations for our results produce some very large and cumbersome representations than the explicit way, which we have followed in this paper.

Now, we begin with some of the fundamental definitions given in Gasper and Rahman[19].

The $q$-shifted factorials are defined in the literature for arbitrary (real or complex) $a, q$ and $b, |q|<1$ as :

$$(a;q)_n = \begin{cases} 1, & n=0 \\ (1\text{-}a)(1\text{-}aq)(1\text{-}aq^2).......(1\text{-}aq^{n-1}), & n \in N \end{cases} \tag{1.2}$$

A generalized basic (or $q$-) hypergeometric series with $r$ numerator parameters $a_1, a_2........a_r$ and $s$ denominator parameters $b_1, b_2........b_s$ □is defined by

$$ {}_r\phi_s(a_1,....a_r;b_1,.....b_s;q,z) \equiv$$

$$ {}_r\phi_s\left[\begin{matrix} a_1,....a_r \\ b_1,.....b_s \end{matrix} ; q, z\right] = \sum_{k=0}^{\infty} \frac{(a_1,....a_r;q)_k}{(q,b_1,....b_s;q)_k} \left[(-1)^k q^{\binom{k}{2}}\right]^{1+s-r} z^k \tag{1.3}$$

where

$$\binom{k}{2} = k(k-1)/2 \tag{1.4}$$

The series $(1.3)$ has the property that if we replace $z$ by $z/a_r$ and then let $a_r \to \infty$ in $(1.3)$, then we obtain a series of the form $(1.3)$ with $r$ is replaced by $r-1$. This is not the case for ${}_r\phi_s$ series defined without the factor $\left[(-1)^k q^{\binom{k}{2}}\right]^{1+s-r}$ in Slater and Bailey[23]. The definition $(1.3)$ is used to handle such limit cases. Also, there is no loss of generality because the Bailey's and Slater's series can be obtained from the $r = s+1$ case of $(1.3)$ by choosing $s$ sufficiently large and setting some of the parameters equal to zero. For further information on basic hypergeometric series and different convergence conditions associated with $(1.3)$, we refer Gasper and Rahman[17].

The following $q$-extension of the (Srivastava-Daoust) generalized Lauricella series in $n$ variables is given in Srivastava and Karlsson[3] as:

$$\phi^{A;B';\dots;B^{(n)}}_{C:D';\dots;D^{(n)}}\begin{bmatrix} x_1 \\ \vdots \\ x_n \end{bmatrix} = \phi^{A;B';\dots;B^{(n)}}_{C:D';\dots;D^{(n)}}\left[\begin{matrix} \left[(a):\theta',\dots\dots\theta^{(n)}\right]:\left[(b'):\phi'\right];\dots\dots;\left[(b^{(n)}):\phi^{(n)}\right] \\ \left[(c):\psi',\dots\dots\psi^{(n)}\right]:\left[(d'):\delta'\right];\dots\dots;\left[(d^{(n)}):\delta^{(n)}\right] \end{matrix}\right.$$

$$\left.;q;\ x_1\dots\dots x_n\right] = \sum_{m_1\dots m_n=0}^{\infty} \Omega(m_1,\dots m_n)\frac{x_1^{m_1}}{m_1!}\cdots\frac{x_n^{m_n}}{m_n!}$$

where

$$\Omega(m_1,\dots m_n) = \frac{\prod_{j=1}^{A}(a_j;q)_{m_1\theta_j'+\dots+m_n\theta_j^{(n)}}\prod_{j=1}^{B'}(b_j';q)_{m_1\phi_j'}\dots\dots\prod_{j=1}^{B^{(n)}}(b_j^{(n)};q)_{m_n\phi_j^{(n)}}}{\prod_{j=1}^{C}(c_j;q)_{m_1\psi_j'+\dots+m_n\psi_j^{(n)}}\prod_{j=1}^{D'}(d_j';q)_{m_1\delta_j'}\dots\dots\prod_{j=1}^{D^{(n)}}(d_j^{(n)};q)_{m_n\delta_j^{(n)}}} \qquad (1.5)$$

where the coefficients and variables are so constrained that the multiple series (1.5) converges.

The $q$-Kampé de Fériet series is given by

$$F^{p:q;k}_{l:m;n}\left[\begin{matrix} a_p : b_q; c_k \\ \alpha_l : \beta_m; \gamma_n \end{matrix}\ ;q,\ x,\ y\right] = \sum_{r,s=0}^{\infty}\frac{\prod_{j=1}^{p}(a_j;q)_{r+s}\prod_{j=1}^{q}(b_j;q)_r\prod_{j=1}^{k}(c_j;q)_s}{\prod_{j=1}^{l}(\alpha_j;q)_{r+s}\prod_{j=1}^{m}(\beta_j;q)_r\prod_{j=1}^{n}(\gamma_j;q)_s}\frac{x^r y^s}{(q;q)_r(q;q)_s}$$

$$\left[(-1)^{r+s}q^{\binom{r+s}{2}}\right]^{l-p}\left[(-1)^r q^{\binom{r}{2}}\right]^{1+m-q}\left[(-1)^s q^{\binom{s}{2}}\right]^{1+n-k} \qquad (1.6)$$

For further detail on notation and convergence conditions, we refer to Srivastava and Karlsson[3].

One of the two theorems (theorem $(2.2)$) on extended Bailey transform proved by Joshi and Vyas [22] is stated below as

If

$$\beta_n = \sum_{r=0}^{n} \alpha_r u_{n-r} v_{n+r} t_{n+2r} w_{pn-r} z_{p'n+r}$$

$$\gamma_n = \sum_{r=n}^{\infty} \delta_r u_{r-n} v_{r+n} t_{2n+r} w_{pr-n} z_{p'r+n} \qquad (1.7)$$

then, subjected to convergence conditions,

$$\sum_{n\geq 0} \alpha_n \gamma_n = \sum_{n\geq 0} \beta_n \delta_n \qquad (1.8)$$

where, $\alpha_r, \delta_r, u_r, v_r, w_r, t_r$ and $z_r$ are any functions of r only and $p$ and $p'$ are any arbitrary integers.

Note that, in many of the papers concerning reductions or transformations of Srivastava-Daoust double hypergeometric series, the parameters $\theta, \varphi, \psi$ and $\delta$'s appearing in equation (1.5) are given some particular constant values. For example, see, Srivastava *et al.*[14]. Use of extended Bailey transform allows us to express these parameters in terms of $p$ that can be assigned any arbitrary integer values. Such results with arbitrary values of these parameters have not appeared previously in the literature. Moreover, it is always possible to derive general reduction formulae involving arbitrary bounded sequence $\{\Omega_n\}$ of complex numbers in place $\delta_n$, provided that the involved series are convergent. Further, the obvious and straightforward generalizations of the results of this paper to reductions or transformations of (m+1) fold series to m-fold series can always be developed after getting the idea of applying *q*-Pfaff-Saalschütz summation theorem used in this paper. This paper is divided into three sections. The second section lists all the new reduction formulae. The section three consists of derivations of the reduction formulae listed in section two. The section fourth deals with the special cases related to each of the reduction formulae.

## *q*-ANALOGUES OF THE REDUCTION FORMULAE FOR SRIVASTAVA-DAOUST FUNCTIONS

In this section, we state following fourteen results as *q*-analogues of the reduction formulae for Srivastava-Daoust functions investigated in[21].

$$\Phi^{D+1:1;1}_{G+1:0;0}\left[\begin{matrix}[d_D:1,1],[z:p+1,p]:[a:1];[\tfrac{k}{za}:1] \\ [g_G:1,1],[k:p+1,p]:-;-\end{matrix}\;;q,\tfrac{kx}{za},x\right] =$$

$$_{D+p^p+(p+1)^{(p+1)}+1}\Phi_{G+p^p+(p+1)^{(p+1)}}\left[\begin{matrix}d_D,\Delta(q;p;z),\Delta(q;p+1;\tfrac{k}{a}),\tfrac{k}{z} \\ g_G,\Delta(q;p+1;k),\Delta(q;p;\tfrac{k}{a})\end{matrix}\;;q,x\right] \qquad (2.1)$$

$$\sum_{n,r\ge 0}\frac{(d_D;q)_{n+r}(v;q)_{r+2n}(z;q)_{(p+1)n+pr}}{(g_G;q)_{n+r}(\frac{qvz}{k};q)_n(k;q)_{(p+1)n+pr}}\frac{(-xq^{-n})^r(-x)^n}{(q;q)_n(q;q)_r}q^{-\binom{n}{2}}=$$

$$\sum_{n,\ge 0}\frac{(d_D;q)_n(v;q)_n(z;q)_{pn}(\frac{k}{v};q)_{pn}(\frac{k}{z};q)_n}{(g_G;q)_n(\frac{k}{v};q)_{(p-1)n}(\frac{qvz}{k};q)_n(k;q)_{(p+1)n}}\frac{(-\frac{vzq}{k}x)^n}{(q;q)_n}q^{\binom{n}{2}} \qquad (2.2)$$

$$\Phi^{D+1:1;1}_{G+1:0;0}\left[\begin{matrix}[d_D:1,1],[w:p-1,p]:[a:1];[\frac{j}{aw}:1];\\ [g_G:1,1],[j:p-1,p]:-;-\end{matrix};q,\frac{x}{a},x\right]=$$

$${}_{D+p^p+(p-1)^{(p-1)}+1}\Phi_{G+p^p+(p-1)^{(p-1)}}\left[\begin{matrix}d_D,\Delta(q;p;wa),\Delta(q;p-1;w),\frac{j}{w}\\ g_G,\Delta(q;p-1;wa),\Delta(q;p;j)\end{matrix};q,\frac{x}{a}\right] \qquad (2.3)$$

$$\sum_{n,r\ge 0}\frac{(d_D;q)_{n+r}(v;q)_{r+2n}(w;q)_{(p-1)n+pr}}{(g_G;q)_{n+r}(\frac{qvw}{j};q)_n(j;q)_{(p-1)n+pr}}\frac{(xq^{-n})^r(-\frac{w}{j}x)^n}{(q;q)_n(q;q)_r}q^{-\binom{n}{2}}=$$

$$\sum_{n\ge 0}\frac{(d_D;q)_n(v;q)_n(wv;q)_{(p+1)n}(w;q)_{(p-1)n}(\frac{j}{w};q)_n}{(g_G;q)_n(wv;q)_{pn}(\frac{qvw}{j};q)_n(j;q)_{pn}}\frac{(-\frac{w}{j}x)^n}{(q;q)_n}q^{-\binom{n}{2}} \qquad (2.4)$$

$$\sum_{n,r\ge 0}\frac{(d_D;q)_{n+r}(z;q)_{(p+2)n+(p+1)r}}{(g_G;q)_{n+r}(h;q)_n(\frac{q^2z}{jh};q)_n(j;q)_{(p-1)n+pr}}\frac{(x)^r(\frac{q}{j}x)^n}{(q;q)_n(q;q)_r}q^{-pn^2-(p+1)nr}=$$

$$\sum_{n\ge 0}\frac{(d_D;q)_n(z;q)_{(p+1)n}(\frac{qz}{h};q)_{(p+1)n}(\frac{jh}{q};q)_{(p+1)n}}{(g_G,h,\frac{q^2z}{jh};q)_n(j,\frac{qz}{h},\frac{jh}{q};q)_{pn}}\frac{(\frac{q}{j}x)^n}{(q;q)_n}q^{-pn^2} \qquad (2.5)$$

$$\sum_{n,r\ge 0}\frac{(d_D;q)_{n+r}(z;q)_{(p+2)n+(p+3)r}}{(g_G;q)_{n+r}(f;q)_{r+2n}(\frac{q^2z}{jv};q)_n(j;q)_{(p-1)n+pr}}\frac{(x)^r(\frac{q}{j}x)^n}{(q;q)_n(q;q)_r}q^{-pn^2-(p+1)nr}=$$

$$\sum_{n\ge 0}\frac{(d_D;q)_n(z;q)_{(p+2)n}(\frac{qz}{f};q)_{(p+1)n}(\frac{jf}{q};q)_{(p+2)n}}{(g_G,\frac{q^2z}{jf};q)_n(j,\frac{qz}{f};q)_{pn}(\frac{jf}{q};q)_{(p+1)n}(f;q)_{2n}}\frac{(\frac{q}{j}x)^n}{(q;q)_n}q^{-pn^2} \qquad (2.6)$$

$$\sum_{n,r\ge 0}\frac{(d_D;q)_{n+r}(z;q)_{(p+2)n+(p+1)r}(\frac{fj}{zq};q)_r}{(g_G;q)_{n+r}(f;q)_{r+2n}(\frac{q^2z}{jv};q)_n(j;q)_{(p-1)n+pr}}\frac{(x)^r(-\frac{fq^{-1/2}}{z}x)^n}{(q;q)_n(q;q)_r}q^{(1/2-p)n^2-pnr}=$$

$$\sum_{n\ge 0}\frac{(d_D;q)_n(z;q)_{(p+1)n}(\frac{qz}{f};q)_{pn}(\frac{jf}{q};q)_{(p+2)n}}{(g_G;q)_n(\frac{qz}{f};q)_{(p-1)n}(\frac{jf}{q};q)_{(p+1)n}(f;q)_{2n}(j;q)_{pn}}\frac{(-\frac{f}{z}x)^n}{(q;q)_n}q^{\binom{n}{2}-pn^2} \qquad (2.7)$$

$$\sum_{n,r\ge 0}\frac{(d_D;q)_{n+r}(z;q)_{(p+4)n+(p+3)r}}{(g_G;q)_{n+r}(f;q)_{r+2n}(\frac{q^2z}{jf};q)_{r+2n}(j;q)_{(p-1)n+pr}}\frac{(x)^r(\frac{q}{j}x)^n}{(q;q)_n(q;q)_r}q^{-pn^2-(p+1)nr}=$$

$$\sum_{n\ge 0}\frac{(d_D;q)_n(z;q)_{(p+3)n}(\frac{qz}{f};q)_{(p+2)n}(\frac{jf}{q};q)_{(p+2)n}}{(g_G;q)_n(f,\frac{q^2z}{jf};q)_{2n}(j;q)_{pn}(\frac{qz}{f};q)_{(p+1)n}(\frac{jf}{q};q)_{(p+1)n}}\frac{(\frac{q}{j}x)^n}{(q;q)_n}q^{-pn^2} \qquad (2.8)$$

$$\sum_{n,r\ge 0}\frac{(d_D;q)_{n+r}(z;q)_{(p+1)n+pr}(\frac{hj}{zq};q)_r}{(g_G;q)_{n+r}(h;q)_n(j;q)_{(p-1)n+pr}}\frac{(x)^r(-\frac{hq^{-1/2}}{z}x)^n}{(q;q)_n(q;q)_r}q^{(1/2-p)n^2-pnr}=$$

$$\sum_{n\ge 0}\frac{(d_D;q)_n(z;q)_{pn}(\frac{qz}{h};q)_{pn}(\frac{jh}{q};q)_{(p+1)n}}{(g_G;q)_n(\frac{qz}{h};q)_{(p-1)n}(\frac{jh}{q};q)_{pn}(h;q)_n(j;q)_{pn}}\frac{(-\frac{h}{z}x)^n}{(q;q)_n}q^{\binom{n}{2}-pn^2} \qquad (2.9)$$

$$\sum_{n,r\ge 0}\frac{(d_D;q)_{n+r}(z;q)_{pn+(p-1)r}(u;q)_r(\frac{j}{zu};q)_r}{(g_G;q)_{n+r}(j;q)_{(p-1)n+pr}}\frac{(x)^r(\frac{x}{z})^n}{(q;q)_n(q;q)_r}q^{(1-p)n^2-(1-p)nr}=$$

$$\sum_{n\ge 0}\frac{(d_D;q)_n(z;q)_{(p-1)n}(uz;q)_{pn}(\frac{j}{u};q)_{pn}}{(g_G;q)_n(uz;q)_{(p-1)n}(\frac{j}{u};q)_{(p-1)n}(j;q)_{pn}}\frac{(x)^n}{(q;q)_n}q^{(1-p)n^2} \qquad (2.10)$$

$$\sum_{n,r\ge 0}\frac{(d_D;q)_{n+r}(w;q)_{(p-1)n+pr}(v;q)_{r+2n}}{(g_G;q)_{n+r}(\frac{qwv}{k};q)_r(k;q)_{(p-1)r+pn}}\frac{(x)^r(-\frac{kq^{-1/2}}{v}x)^n}{(q;q)_n(q;q)_r}q^{(p-3/2)n^2-(p-2)nr}=$$

$$\sum_{n\ge 0}\frac{(d_D;q)_n(v;q)_n(\frac{k}{v};q)_{(p-1)n}(w;q)_{(p-1)n}(wv;q)_{(p+1)n}}{(g_G;q)_n(\frac{qwv}{k};q)_n(wv;q)_{pn}(\frac{k}{v};q)_{(p-2)n}(k;q)_{pn}}\frac{(x)^n}{(q;q)_n} \qquad (2.11)$$

$$\sum_{n,r\ge 0}\frac{(d_D;q)_{n+r}(w;q)_{(p-1)n+pr}(v;q)_{r+2n}(\frac{k}{vw};q)_n}{(g_G;q)_{n+r}(k;q)_{(p+1)n+pr}}\frac{(x)^r(xw)^n}{(q;q)_n(q;q)_r}q^{(p-1)(n^2+nr)}=$$

$$\sum_{n\ge 0}\frac{(d_D;q)_n(v;q)_n(\frac{k}{v};q)_{pn}(w;q)_{(p-1)n}(wv;q)_{(p+1)n}}{(g_G;q)_n(wv;q)_{pn}(\frac{k}{v};q)_{(p-1)n}(k;q)_{(p+1)n}}\frac{(x)^n}{(q;q)_n} \qquad (2.12)$$

$$\sum_{n,r\ge 0}\frac{(d_D;q)_{n+r}(w;q)_{(p-1)n+pr}(\frac{ke}{wq};q)_n}{(g_G;q)_{n+r}(e;q)_r(k;q)_{(p-2)r+(p-1)n}}\frac{(x)^r(-\frac{wq^{1/2}}{e}x)^n}{(q;q)_n(q;q)_r}q^{(p-3/2)n^2-(p-2)nr}=$$

$$\sum_{n\ge 0}\frac{(d_D;q)_n(\frac{ke}{q};q)_{pn}(w;q)_{(p-1)n}(\frac{wq}{e};q)_{(p-1)n}}{(g_G;q)_n(e;q)_n(\frac{ke}{q};q)_{(p-1)n}(\frac{wq}{e};q)_{(p-2)n}(k;q)_{(p-1)n}}\frac{(x)^n}{(q;q)_n} \qquad (2.13)$$

$$\sum_{n,r\ge 0}\frac{(d_D;q)_{n+r}(a^3;q)_{r+3n}}{(g_G;q)_{n+r}(h;q)_n(\frac{a^6q^3}{h^2};q^2)_n}\frac{(-xq)^n x^r}{(q;q)_n(q^2;q^2)_r}q^{-\frac{3}{2}n^2-\frac{r^2}{2}-3nr}=$$

$$\sum_{n\ge 0}\frac{(d_D;q)_n(a^3;q)_n(\frac{h^2}{a^3q};q)_n(\frac{q^2a^3}{h^2};q)_n}{(g_G;q)_n(h^2;q^2)_n(\frac{q^3a^6}{h^2};q^2)_n}\frac{(a^3x)^n}{(q^2;q^2)_n}q^{n/2} \qquad (2.14)$$

## DERIVATIONS OF THE RESULTS FROM (2.1) TO (2.14)

To derive the different results mentioned in previous section, we decide different expressions for $\alpha_r, \delta_r, u_r, v_r, w_r, t_r$ and $z_r$ in extended Bailey transform, which yields closed form for $\beta_n$ and $\delta_n$ by utilizing $q$-Pfaff-Saalschütz summation theorem as mentioned in Gasper and Rahman [17]. The final results are obtained with help of equation $(1.8)$.

We follow the same aforementioned process to obtain each of the reduction formula listed in section two. Note that, $D$ and $G$ are positive integers, while $p$ and $p'$ are arbitrary integers.

(i). Choose $\alpha_n = \dfrac{(a;q)_n (\frac{kx}{az})^n}{(q;q)_n}$, $u_n = \dfrac{(\frac{k}{az};q)_n x^n}{(q;q)_n}$, $z_n = \dfrac{(z;q)_n}{(k;q)_n}$, and $p = p'$.

(ii). Select $\alpha_n = \dfrac{(-xq^{1/2})^n}{(\frac{qvz}{k};q)_n (q;q)_n}$, $u_n = \dfrac{x^n}{(q;q)_n}$, $z_n = \dfrac{(z;q)_n}{(k;q)_n}$, $v_n = (v;q)_n$ and $p = p'$

(iii). Select $\alpha_n = \dfrac{(a;q)_n (\frac{x}{a})^n}{(q;q)_n}$, $u_n = \dfrac{(\frac{j}{aw};q)_n x^n}{(q;q)_n}$, $w_n = \dfrac{(w;q)_n}{(j;q)_n}$ and $p = p'$

(iv). Select $\alpha_n = \dfrac{(-\frac{w}{j} xq^{1/2})^n}{(\frac{qvw}{j};q)_n (q;q)_n}$, $u_n = \dfrac{x^n}{(q;q)_n}$, $w_n = \dfrac{(w;q)_n}{(j;q)_n}$, $v_n = (v;q)_n$ and $p = p'$

(v). Select

$$\alpha_n = \frac{(\frac{q}{j} x)^n}{(\frac{q^2 z}{jh};q)_n (h;q)_n (q;q)_n}, \quad u_n = \frac{x^n}{(q;q)_n}, \quad w_n = \frac{1}{(j;q)_n}, \quad z_n = (z;q)_n \text{ and } p' = p+1$$

(vi). Select

$$\alpha_n = \frac{(\frac{q}{j} x)^n}{(\frac{q^2 z}{jv};q)_n (q;q)_n}, \quad u_n = \frac{x^n}{(q;q)_n}, \quad w_n = \frac{1}{(j;q)_n},$$

$$v_n = \frac{1}{(f;q)_n}, \quad z_n = (z;q)_n, \text{ and } p' = p+2$$

(vii). Select

$$\alpha_n = \frac{(-\frac{fq^{-1/2}}{z} x)^n}{(q;q)_n}, \quad u_n = \frac{(\frac{fj}{zq};q)_n x^n}{(q;q)_n}, \quad w_n = \frac{1}{(j;q)_n},$$

$$v_n = \frac{1}{(f;q)_n}, \quad z_n = (z;q)_n, \text{ and } p' = p+1$$

(viii). Select

$$\alpha_n = \frac{(\frac{q}{j} x)^n}{(q;q)_n}, \quad u_n = \frac{x^n}{(q;q)_n}, \quad w_n = \frac{1}{(j;q)_n},$$

$$v_n = \frac{1}{(\frac{q^2 z}{fj};q)_n (f;q)_n}, \quad z_n = (z;q)_n, \text{ and } p' = p+3$$

(ix). Select $\alpha_n = \dfrac{(-\frac{hq^{-1/2}}{z} x)^n}{(h;q)_n (q;q)_n}$, $u_n = \dfrac{(\frac{hj}{zq};q)_n x^n}{(q;q)_n}$, $w_n = \dfrac{1}{(j;q)_n}$, $z_n = (z;q)_n$, and $p' = p$

(x). Choose

$$\alpha_n=\frac{(\frac{x}{z})^n}{(q;q)_n},\ u_n=\frac{(u;q)_n(\frac{j}{zu};q)_n x^n}{(q;q)_n},\quad w_n=\frac{1}{(j;q)_n}, z_n=(z;q)_n,\ \text{and}\ p'=p-1$$

$$\alpha_n=\frac{(-\frac{kq^{-1/2}}{v}x)^n}{(q;q)_n},\ u_n=\frac{x^n}{(\frac{qwv}{k};q)_n(q;q)_n},\quad w_n=(w;q)_n,$$

(xi). Select

$$z_n=\frac{1}{(k;q)_n},\ v_n=(v;q)_n,\ \text{and}\ p'=p-1$$

(xii). Choose

$$\alpha_n=\frac{(\frac{k}{vw};q)_n(wx)^n}{(q;q)_n},\ u_n=\frac{x^n}{(q;q)_n},\ v_n=(v;q)_n, w_n=(w;q)_n,$$

$$z_n=\frac{1}{(k;q)_n},\ \text{and}\ p'=p$$

(xiii). Consider

$$\alpha_n=\frac{(\frac{ke}{qw};q)_n(-\frac{wq^{1/2}}{e}x)^n}{(q;q)_n},\ u_n=\frac{x^n}{(e;q)_n(q;q)_n},\quad w_n=(w;q)_n,$$

$$z_n=\frac{1}{(k;q)_n},\ \text{and}\ p'=p-2$$

(xiv). Choose

$$\alpha_n=\frac{(-xq)^n}{(h^2;q^2)_n\left(\frac{a^6q^3}{h^2};q^2\right)_n(q;q)_n},\ u_n=\frac{x^n}{(q^2;q^2)_n},\quad t_n=(a^3;q)_n$$

## PARTICULAR CASES OF DERIVED REDUCTION FORMULAE

By assigning different integer values to arbitrary variable $p, D$ and $G$, we obtain the $q$-analogues of well-known results like Kampé de Fériet reduction formula, Euler transformation formula and Whipple quadratic transformation formula (which is also known as Sears-Carlitz transformation formula) as recorded in Srivastava and Karlsson[3], Andrews *et al.*[15] and Gasper and Rahman[17] respectively. When $q$ tends to 1, the aforementioned well-known results transform into ordinary Kampé de Fériet reduction formula discussed in Srivastava and Karlsson[3], Euler transformation formula[15] and Whipple quadratic transformation formula[15] respectively (see[24,25] also).

(i). In equation (2.1) selecting $p=0$, we obtain

$$\sum_{n,r\geq 0}\frac{(d_D;q)_{n+r}(a;q)_n(\frac{k}{az};q)_r(z;q)_n}{(g_G;q)_{n+r}(k;q)_n}\frac{(\frac{kx}{az})^n(x)^r}{(q;q)_n(q;q)_r}=$$

$$\sum_{n,\geq 0}\frac{(d_D;q)_n(\frac{k}{z};q)_n(\frac{k}{a};q)_n}{(g_G;q)_n(k;q)_n}\frac{(x)^n}{(q;q)_n}\qquad(4.1)$$

As $q$ tends to 1, the above result converts into the ordinary Kampé de Fériet reduction formula[3] .

(ii). For the choice $p = 0 = D = G$ in equation $(2.2)$, we obtain Sears-Carlitz transformation formula from ${}_3\varphi_2$ to ${}_5\varphi_4$ as given below.

$$
{}_3\varphi_2\left[\begin{matrix} v, \frac{k}{z}, \frac{qv}{k} \\ k, \frac{qvz}{k} \end{matrix} ; q, \ xz\right] =
$$

$$
\frac{(vx;q)_\infty}{(x;q)_\infty}\, {}_5\varphi_4\left[\begin{matrix} \sqrt{v}, -\sqrt{v}, \sqrt{vq}, -\sqrt{vq}, z \\ \frac{qvz}{k}, k, vx, \frac{q}{x} \end{matrix} ; q, \ q\right] \qquad (4.2)
$$

Taking $q$ tends to 1; we obtain ordinary Whipple's transformation formula[25].

(iii). For $p = 1$ in equation $(2.3)$, we again obtain a $q$-analogue of the ordinary Kampé de Fériet reduction formula[3] and the choice $p = 0 = D = G$ gives $q$-analogue of Euler transformation formula[25].

(iv). When $p = 0 = D = G$ in equation $(2.4)$ we obtain Sears-Carlitz transformation formula like equation (4.2) and $p = 1$, $D = G = 0$ gives a $q$-analogue of a reduction formula for Horn's $H_4$ hypergeometric series[3] as

$$
\sum_{n,r\geq 0} \frac{(v;q)_{2n+r}\,(w;q)_r}{(\frac{qwv}{j};q)_n\,(j;q)_r} \frac{(xq^{-n})^r\,(-\frac{xw}{j})^n}{(q;q)_n\,(q;q)_r}\, q^{-\binom{n}{2}} =
$$

$$
\sum_{n,\geq 0} \frac{(v;q)_n\,(wv;q)_{2n}\,(\frac{j}{w};q)_n}{(j;q)_n\,(wv;q)_n\,(\frac{qwv}{j};q)_n} \frac{(-\frac{w}{j}x)^n}{(q;q)_n}\, q^{\binom{n}{2}} \qquad (4.3)
$$

As $q$ tends to 1, we recover the reduction formula for Horn's $H_4$ hypergeometric series[3].

(v). When $p = 0 = D = G$ in equation $(2.5)$, we again obtain Sears-Carlitz transformation[17] formula given by equation $(4.2)$. As $q$ tends to 1, we obtain Whipple quadratic formula [15].

(vi). In equation $(2.9)$ selecting $p = 0$, we get a $q$-analogue of a Kampé de Fériet reduction formula[3].

(vii). For $p = 0$ or $p = 1$ in equation $(2.10)$, the $q$-analogue of a Kampé de Fériet reduction formula[3] follows.

(viii). For $p = 0 = D = G$ in equation $(2.12)$, we obtain Sears-Carlitz transformation formula[17] as given in equation $(4.2)$, which leads to Whipple quadratic formula[15], as $q$ tends to 1.

*** Corresponding author: *Yashoverdhan Vyas***

Department of Mathematics, School of Engineering,
Sir Padampat Singhania University, Udaipur, Rajasthan.

*Email: yashoverdhan.vyas@spsu.ac.in*